\documentclass[a4paper,11pt,twoside]{article}

\RequirePackage[a4paper,head=1cm]{geometry}
\usepackage[utf8x]{inputenc}
\usepackage{amsmath}
\usepackage{amssymb}
\usepackage{hyperref}
\usepackage{dsfont}
\newtheorem{Definition}{Definition}[section]
\newtheorem{Theoreme}{Theorem}

\newtheorem{Proposition}{Proposition}[section]
\newtheorem{Remarque}{\bf Remark}

\title{\bf A remark on Besov spaces interpolation over the 2-adic group} 
\author{Diego Chamorro } 

\begin{document} 
\maketitle 
\begin{center}
\begin{minipage}[l]{150mm}
\begin{scriptsize}\abstract{Motivated by a recent result which identifies in the special setting of the $2$-adic group the Besov space $\dot{B}^{1,\infty}_{1}(\mathbb{Z}_2)$ with $BV(\mathbb{Z}_2)$, the space of function of bounded variation, we study in this article some functional relationships between Besov spaces.\\[3mm]
}\end{scriptsize}
\end{minipage}
\end{center}
\section{Introduction}

The starting point of this article is given by the following inequality proved by A. Cohen, W. Dahmen, I. Daubechies \& R. De Vore in \cite{Cohen2}.
For a function $f:\mathbb{R}^n\longrightarrow \mathbb{R}$ such that $f\in BV \cap \dot{B}^{-1,\infty}_{\infty}$ we have
\begin{equation}\label{EQ1}
\|f\|_{L^2}^2\leq C\|f\|_{BV}\|f\|_{\dot{B}^{-1,\infty}_{\infty}}
\end{equation}
Here $BV$ denotes the space of functions of bounded variation and $\dot{B}^{-1,\infty}_{\infty}$ stands for an homogeneous Besov space.
In the article \cite{Chame2}, we proved that in the special setting of the $2$-adic group $\mathbb{Z}_2$, the space $BV(\mathbb{Z}_2)$ can be identified to the Besov space $\dot{B}^{1,\infty}_{1}(\mathbb{Z}_2)$ and therefore, inequality (\ref{EQ1}) becomes
\begin{equation}\label{EQ2}
\|f\|_{L^2}^2\leq C\|f\|_{\dot{B}^{1,\infty}_{1}}\|f\|_{\dot{B}^{-1,\infty}_{\infty}}
\end{equation}
Note that the previous estimate is \textit{false} in $\mathbb{Z}_2$, see \cite{Chame2} for a counterexample. The identification between these two functional spaces and the consequences on the inequality (\ref{EQ1}) are very surprising in the sense that these estimates \textit{depend} on the underlying group structure: compare the topological properties of $\mathbb{R}^n$ to the totally discontinuous setting of $\mathbb{Z}_2$.\\

However, one may think that the Besov norm $\|\cdot\|_{\dot{B}^{1,\infty}_{1}}$ in the right hand side of (\ref{EQ2}) is \textit{too} small to achieve the inequality. Thus, it is a natural question to study the validity of (\ref{EQ2}) if we replace this norm by a \textit{bigger} one (just think on the inclusion of Besov spaces $\dot{B}^{1,q}_{1}\subset\dot{B}^{1,\infty}_{1}$ valid for $q\geq 1$).
The answer to this question is given by the next result

\begin{Theoreme}\label{Theo1}
If $f:\mathbb{Z}_2\longrightarrow \mathbb{R}$ is a function such that $f\in\dot{B}^{1,q}_{1}\cap \dot{B}^{-1,\infty}_{\infty}(\mathbb{Z}_2)$ with $q>2$, then the following inequality is false:
\begin{equation}\label{EQ3}
\|f\|_{L^2}^2\leq C\|f\|_{\dot{B}^{1,q}_{1}}\|f\|_{\dot{B}^{-1,\infty}_{\infty}}
\end{equation}
\end{Theoreme}
This is the main theorem of this article and we will construct a counterexample in the section \ref{ProofTH1} below, but before, it would be interesting to compare inequality (\ref{EQ3}) to the general estimates given by the interpolation theory\footnote{see the book \cite{Bergh} for more details.}.\\ 

Indeed, following this general theory, we can obtain inequalities of the form
\begin{equation}\label{EQ4}
\|f\|_{L^2}^2\leq C\|f\|_{\dot{B}^{s_0,q_0}_{p_0}}\|f\|_{\dot{B}^{s_1,q_1}_{p_1}}
\end{equation}
for some special values of the real parameters $s_0,s_1,p_0,p_1,q_0,q_1$.\\ 

Perhaps the most popular case is given by the real method: set $p_0=p_1=p$, fix $0<\theta<1$ and suppose $s_0\neq s_1$ with the relationship $s=(1-\theta)s_0+ \theta s_1$. We obtain the following expression
$$\big(\dot{B}^{s_0,q_0}_{p},\dot{B}^{s_1,q_1}_{p}\big)_{\theta,q}=\dot{B}^{s,q}_{p}$$
which gives us the estimate
\begin{equation}\label{EQ5}
\|f\|_{\dot{B}^{s,q}_{p}}\leq C\|f\|_{\dot{B}^{s_0,q_0}_{p}}^{1-\theta}\|f\|_{\dot{B}^{s_1,q_1}_{p}}^{\theta}
\end{equation}
It is very important to remark that in this particular case no relationship between $q_0,q_1$ and $q$ is asked. 
Obviously, inequality (\ref{EQ3}) can not be obtained from (\ref{EQ5}), since $p_0 \neq p_1$.\\

The case when $p_0\neq p_1$ is more restrictive and following the complex method we have for $1\leq p_0,q_0\leq +\infty$ and $1\leq p_1,q_1<+\infty$ the formula
$$\big[\dot{B}^{s_0,q_0}_{p_0},\dot{B}^{s_1,q_1}_{p_1}\big]_{\theta}=\dot{B}^{s,q}_{p}$$
which gives us an estimate of the type (\ref{EQ4}) with $s=(1-\theta)s_0+ \theta s_1$, $\frac{1}{p}=\frac{1-\theta}{p_0}+\frac{\theta}{p_1}$, and $\frac{1}{q}=\frac{1-\theta}{q_0}+\frac{\theta}{q_1}$. Note that we have in this case a relationship between $q_0,q_1$ and $q$. Again, this method can not be applied to inequality (\ref{EQ3}). \\

It seems of course that inequality (\ref{EQ3}) cannot be obtained by an simple interplation argument -actually this inequality is false in $\mathbb{R}^n$-, but what it would make it plausible in the setting of $\mathbb{Z}_2$ is the special relationship between inequalities (\ref{EQ1}) and (\ref{EQ2}) and this is the main reason why theorem \ref{Theo1} is relevant.\\

The plan of the article is the following. In section \ref{IntroP} we recall some properties of the $p$-adic spaces, in section \ref{ESPFUNC} we give the definition of Besov spaces over the $2$-adic group $\mathbb{Z}_2$ and in section \ref{ProofTH1} we prove theorem \ref{Theo1}.

\section{$p$-adic groups}\label{IntroP}
Our main reference here are the books \cite{VVZ}, \cite{Koblitz} and \cite{Amice} where more details concerning the topological structure of the $p$-adic groups can be found.\\

We write $a|b$ when $a$ divide $b$ or, equivalently, when $b$ is a multiple of $a$. Let $p$ be any prime number, for $0\neq x\in \mathbb{Z}$, we define the $p$-adic valuation of $x$ by $\gamma(x)=\max\{r: p^{r}|x\}\geq 0$ and, for any rational number $x=\frac{a}{b}\in \mathbb{Q}$, we write $\gamma(x)=\gamma(a)-\gamma(b)$. Furthermore if $x=0$, we agree to write $\gamma(0)=+\infty$.\\ 

Let $x \in \mathbb{Q}$ and $p$ be any prime number, with the $p$-adic valuation of $x$ we can construct a norm by writing
\begin{equation}\label{Normeadique}
|x|_{p}=\left\lbrace \begin{array}{ll}
p^{-\gamma} &\mbox{if}\quad x\neq 0\\[5mm]
p^{-\infty}=0 &\mbox{if}\quad x= 0.
\end{array}\right.
\end{equation}
This expression satisfy the following properties
\begin{enumerate}
\item[a)]$|x|_{p}\geq 0$, and $|x|_{p}=0 \iff x=0$;
\item[b)] $|xy|_{p}=|x|_{p}|y|_{p}$;
\item[c)] $|x+y|_{p}\leq \max\{|x|_{p},|y|_{p}\}$, with equality when $|x|_{p}\neq |y|_{p}$.
\end{enumerate}
When a norm satisfy $c)$ it is called a non-Archimedean norm and an interesting fact is that over $\mathbb{Q}$ \textit{all} the possible norms are equivalent to $|\cdot|_{p}$ for some $p$: this is the so-called Ostrowski theorem, see \cite{Amice} for a proof. 
\begin{Definition}
Let $p$ be a any prime number. We define the field of $p$-adic numbers $\mathbb{Q}_{p}$ as the completion of $\mathbb{Q}$ when using the norm $|\cdot|_{p}$.
\end{Definition}
We present in the following lines the algebraic structure of the set $\mathbb{Q}_{p}$. Every $p$-adic number $x\neq 0$ can be represented in a unique manner by the formula
\begin{equation}\label{Hensel}
x=p^{\gamma}(x_{0}+x_{1}p+x_{2}p^{2}+...),
\end{equation}
where $\gamma=\gamma(x)$ is the $p$-adic valuation of $x$ and $x_{j}$ are integers such that $x_{0}>0$ and $0\leq x_{j}\leq p-1$ for $j=1,2,...$. Remark that this canonical representation implies the identity $|x|_{p}=p^{-\gamma}$.\\

Let $x,y \in \mathbb{Q}_{p}$, using the formula (\ref{Hensel}) we define the sum of $x$ and $y$ by 
$x+y=p^{\gamma (x+y)}(c_{0}+c_{1}p+c_{2}p^{2}+...)$  with $0\leq c_{j}\leq p-1$ and $c_{0}>0$, where $\gamma(x+y)$ and $c_{j}$ are the unique solution of the equation
$$p^{\gamma(x)}(x_{0}+x_{1}p+x_{2}p^{2}+...)+ p^{\gamma(y)}(y_{0}+y_{1}p+y_{2}p^{2}+...)=p^{\gamma(x+y)}(c_{0}+c_{1}p+c_{2}p^{2}+...).$$
Furthermore, for $a,x \in \mathbb{Q}_{p}$, the equation $a+x=0$ has a unique solution in $\mathbb{Q}_{p}$ given by $x=-a$. In the same way, the equation $ax=1$ has a unique solution in $\mathbb{Q}_{p}$: $x=1/a$.\\

We take now a closer look at the topological structure of $\mathbb{Q}_{p}$. With the norm $|\cdot|_p$ we construct a distance over $\mathbb{Q}_p$ by writing
\begin{equation}\label{dist}
d(x,y)=|x-y|_{p} 
\end{equation} 
and we define the balls $B_{\gamma}(x)=\left\{y\in \mathbb{Q}_{p}: \; d(x,y)\leq p^{\gamma}\right\}$  with $\gamma \in \mathbb{Z}$. Remark that, from the properties of the $p$-adic valuation, this distance has the \textit{ultra-metric} property (\textit{i.e.} $d(x,y)\leq \max\{d(x,z),d(z,y)\}\leq |x|_{p}+|y|_{p}$).\\

We gather with the next proposition some important facts concerning the balls in $\mathbb{Q}_{p}$.
\begin{Proposition}
Let $\gamma$ be an integer, then we have
\begin{enumerate}
\item[1)] the ball $B_{\gamma}(x)$ is a open and a closed set for the distance (\ref{dist}).
\item[2)] every point of $B_{\gamma}(x)$ is its center.
\item[3)] $\mathbb{Q}_{p}$ endowed with this distance is a complete Hausdorff metric space.
\item[4)] $\mathbb{Q}_{p}$ is a locally compact set.
\item[5)] the $p$-adic group $\mathbb{Q}_{p}$ is a totally discontinuous space.
\end{enumerate}
\end{Proposition}
\section{Functional spaces}\label{ESPFUNC}
In this article, we will work with the subset $\mathbb{Z}_{2}$ of $\mathbb{Q}_{2}$ which is defined by $\mathbb{Z}_{2}=\{x\in \mathbb{Q}_{2}:\; |x|_{2}\leq 1\}$, and we will focus on real-valued functions over $\mathbb{Z}_{2}$. Since $\mathbb{Z}_{2}$ is a locally compact commutative group, there exists a Haar measure $dx$ which is translation invariant $i.e.$: $d(x+a)=dx$, furthermore we have the identity $d(xa)=|a|_{2}dx$ for $a\in \mathbb{Z}_{2}^{*}$. We will normalize the measure $dx$ by setting
$$\int_{\{|x|_{2}\leq1\}}dx=1.$$
This measure is then unique and we will note $|E|$ the measure for any subset $E$ of $\mathbb{Z}_{2}$. \\

Lebesgue spaces $L^{p}(\mathbb{Z}_{2})$ are thus defined in a natural way: $\|f\|_{L^p}=\left(\displaystyle{\int_{\mathbb{Z}_{2}}}|f(x)|^{p}dx\right)^{1/p}$ for $1\leq p< +\infty$, with the usual modifications when $p=+\infty$.\\

Let us now introduce the Littlewood-Paley decomposition in $\mathbb{Z}_{2}$. We note $\mathcal{F}_{j}$ the Boole algebra formed by the equivalence classes $E\subset \mathbb{Z}_{2}$ modulo the sub-group $2^{j}\mathbb{Z}_{2}$. Then, for any function $f\in L^{1}(\mathbb{Z}_{2})$, we call $S_{j}(f)$ the conditionnal expectation of $f$ with respect to $\mathcal{F}_{j}$:
$$S_{j}(f)(x)=\frac{1}{|B_{j}(x)|}\int_{B_{j}(x)}f(y)dy.$$
The dyadic blocks are thus defined by the formula $\Delta_{j}(f)=S_{j+1}(f)-S_{j}(f)$ and the Littlewood-Paley decomposition of a function $f:\mathbb{Z}_{2}\longrightarrow \mathbb{R}$ is given by
\begin{equation}\label{LPP}
f=S_{0}(f)+\sum_{j=0}^{+\infty}\Delta_{j}(f)\qquad \mbox{where } S_{0}(f)=\int_{\mathbb{Z}_{2}}f(x)dx.
\end{equation}

We will need in the sequel some very special sets noted $Q_{j,k}$. Here is the definition and some properties:
\begin{Proposition}\label{pala2bis}
Let $j\in \mathbb{N}$ and $k=\{0,1,...,2^{j}-1\}$. Define the subset $Q_{j,k}$ of $\mathbb{Z}_2$ by
\begin{equation}\label{QJKset}
Q_{j,k}=\left\{k+2^{j}\mathbb{Z}_{2}\right\}.
\end{equation}
Then
\begin{enumerate}
\item[1)] We have the identity $\mathcal{F}_{j}=\underset{0\leq k<2^{j}}{\bigcup} Q_{j,k}$,
\item[2)] For $k=\{0,1,...,2^{j}-1\}$ the sets $Q_{j,k}$ are mutually disjoint, 
\item[3)] $|Q_{j,k}|=2^{-j}$ for all $k$,
\item[4)] the $2$-adic valuation is constant over $Q_{j,k}$.
\end{enumerate}
\end{Proposition}

The verifications are easy and left to the reader.\\

With the Littlewood-Paley decomposition given in (\ref{LPP}), we obtain the following equivalence for the Lebesgue spaces $L^p(\mathbb{Z}_2)$ with $1<p<+\infty$:
$$\|f\|_{L^p}\simeq \|S_0(f)\|_{L^p}+\left\|\bigg(\sum_{j\in \mathbb{N}}|\Delta_{j}f|^{2}\bigg)^{1/2} \right\|_{L^p}.$$
See the book \cite{Stein0}, chapter IV, for a general proof. \\

For Besov spaces we will define them by the norm
\begin{equation}\label{besov1}
\|f\|_{B^{s,q}_{p}}\simeq\|S_{0}f\|_{L^p}+\left(\sum_{j\in \mathbb{N}}2^{jsq}\|\Delta_{j}f\|^{q}_{L^p}\right)^{1/q}
\end{equation}
where $s\in \mathbb{R}$, $1\leq p,q<+\infty$ with the necessary modifications when $p,q=+\infty$.\\
\begin{Remarque}
\emph{For homogeneous functional spaces $\dot{B}^{s,q}_{p}$, we drop out the term $\|S_{0}f\|_{L^p}$ in (\ref{besov1}).}
\end{Remarque}
\section{Proof of the theorem \ref{Theo1}}\label{ProofTH1}

To begin the construction of the counterexample we consider $0<j_0<j_1$ two integers and we fix $\alpha,\beta\in\mathbb{R}$ such that
\begin{equation}\label{EQ6}
2^{2j_0}\leq \frac{\beta}{\alpha}.
\end{equation}
Take now a decreasing sequence $(\varepsilon_j)_{j\in \mathbb{N}}\in \ell^q(\mathbb{N})$ with $q>2$ such that $\varepsilon_0=1$ and $(\varepsilon_j)_{j\in \mathbb{N}}\notin \ell^2(\mathbb{N})$.

Define $N_j$ in the following form
\begin{equation}\label{DefinitionN}
N_j=
\begin{cases}
2^{j} &\text{if}\quad  0<j<j_0,\\[5mm]
2^{-j}\frac{\beta}{\alpha} &\text{if}\quad  j_0\leq j\leq j_1.
\end{cases}
\end{equation}

We construct a function $f:\mathbb{Z}_2\longrightarrow \mathbb{R}$ by considering his values over the dyadic blocs and we will use for this the sets $Q_{j,k}$ defined in (\ref{QJKset}):
\begin{equation*}
\Delta_{j}f(x)=
\begin{cases}
\varepsilon_j\alpha 2^{j} &\text{over}\quad  Q_{j+1,0},\\
-\varepsilon_j\alpha 2^{j} &\text{over}\quad  Q_{j+1,1},\\
\varepsilon_j\alpha 2^{j} &\text{over}\quad  Q_{j+1,2},\\
-\varepsilon_j\alpha 2^{j} &\text{over}\quad  Q_{j+1,3},\\
&\vdots\\
\varepsilon_j\alpha 2^{j} &\text{over} \quad Q_{j+1,2N_j-2},\\
-\varepsilon_j\alpha 2^{j} &\text{over}\quad  Q_{j+1,2N_j-1},\\
0 &\text{elsewhere}.
\end{cases}
\end{equation*}
Remark that, with this definition of $\Delta_jf(x)$ we have the identities
\begin{itemize}
\item $\|\Delta_jf\|_{L^\infty}=\varepsilon_j\alpha 2^j$,
\item $\|\Delta_jf\|_{L^1}=\varepsilon_j \alpha N_j$,
\item $\|\Delta_jf\|_{L^2}^2=\varepsilon_j^2 \alpha^2 2^jN_j$.
\end{itemize}
From this quantities we construct the following norms
\begin{itemize}
\item[(a)] for the Besov space $\dot{B}^{-1,\infty}_{\infty}$ we have\\

$\|f\|_{\dot{B}^{-1,\infty}_{\infty}}=\underset{j\in \mathbb{N}}{\sup}\;2^{-j}\|\Delta_jf\|_{L^\infty}=\alpha$, since the sequence $(\varepsilon_j)_{j\in\mathbb{N}}$ is decreasing and $\varepsilon_0=1$.
\item[(b)] for the Besov space $\dot{B}^{1,q}_{1}$ we write\\ 

$\|f\|_{\dot{B}^{1,q}_{1}}^q=\displaystyle{\sum_{j=0}^{j_1}}\left(2^j\|\Delta_jf\|_{L^1}\right)^q=\displaystyle{\sum_{j=0}^{j_1}}2^{jq}\varepsilon_j^q\alpha^q N_j^q=\alpha^q\left(\displaystyle{\sum_{j=0}^{j_0}}2^{jq}\varepsilon_j^qN_j^q+\displaystyle{\sum_{j>j_0}^{j_1}}2^{jq}\varepsilon_j^q N_j^q\right)$\\
 
We use now the values of $N_j$ given in (\ref{DefinitionN}) and the relationship (\ref{EQ6}) to obtain
 $$=\alpha^q\left(\displaystyle{\sum_{j=0}^{j_0}}2^{2jq}\varepsilon_j^q+\displaystyle{\sum_{j>j_0}^{j_1}}\varepsilon_j^q\frac{\beta^q}{\alpha^q}\right)=\beta^q\left(\displaystyle{\sum_{j=0}^{j_0}}2^{2jq}\frac{\alpha^q}{\beta^q}\varepsilon_j^q+\displaystyle{\sum_{j>j_0}^{j_1}}\varepsilon_j^q\right)
\simeq \beta^q\left(\sum_{j=0}^{j_0}2^{q(2j-2j_0)}\varepsilon_j^q+\displaystyle{\sum_{j>j_0}^{j_1}}\varepsilon_j^q\right).$$
Then we have $\|f\|_{\dot{B}^{1,q}_{1}}\simeq \beta\left(C_1+\displaystyle{\sum_{j>j_0}^{j_1}}\varepsilon_j^q\right)^{1/q}.$

\item[(c)] For the Lebesgue space $L^2$ we use the same arguments above to obtain\\

$\|f\|_{L^2}^2=\displaystyle{\sum_{j=0}^{j_1}}\varepsilon_j^2\alpha^22^jN_j=\alpha^2\left(\displaystyle{\sum_{j=0}^{j_0}}2^{2j}
\varepsilon_j^2+\displaystyle{\sum_{j>j_0}^{j_1}}\varepsilon_j^2\frac{\beta}{\alpha}\right)\simeq \alpha\beta\left(C_2+\displaystyle{\sum_{j>j_0}^{j_1}}\varepsilon_j^2\right)$.
\end{itemize}
Once these norms are computed, we go back to the inequality 
$$\|f\|_{L^2}^2\leq C\|f\|_{\dot{B}^{1,q}_{1}}\|f\|_{\dot{B}^{-1,\infty}_{\infty}}$$
and we have
$$\alpha\beta\left(C_2+\displaystyle{\sum_{j>j_0}^{j_1}}\varepsilon_j^2\right)\leq C\times\alpha \times \beta\left(C_1+\displaystyle{\sum_{j>j_0}^{j_1}}\varepsilon_j^q\right)^{1/q}.$$
But, by hypothesis, we have $(\varepsilon_j)_{j\in \mathbb{N}}\notin \ell^2(\mathbb{N})$ and $(\varepsilon_j)_{j\in \mathbb{N}}\in \ell^q(\mathbb{N})$, thus, for $j_1$ big enough it is impossible to find an universal constant $C$ such that the above inequality is true. 
\begin{flushright}{$\blacksquare$}\end{flushright}  


\begin{flushright}
\begin{minipage}[r]{80mm}
Diego \textsc{Chamorro}\\[5mm]
Laboratoire d'Analyse et de Probabilités\\ 
Université d'Evry Val d'Essonne \& ENSIIE\\[2mm]
1 square de la résistance,\\
91025 Evry Cedex\\[2mm]
diego.chamorro@m4x.org
\end{minipage}
\end{flushright}

\end{document}